\input amstex
\documentstyle{amsppt}
%%%\magnification=\magstep1
\hsize = 6.5 truein
\vsize = 9 truein
\NoBlackBoxes
\TagsAsMath
 
%Section macros
\newskip\sectionskipamount
\sectionskipamount = 24pt plus 8pt minus 8pt
\def\sectionskip{\vskip\sectionskipamount}
\define\sectionbreak{%
        \par  \ifdim\lastskip<\sectionskipamount
        \removelastskip  \penalty-2000  \sectionskip  \fi}
\define\section#1{%
        \sectionbreak                   %Encourages a page break, else inserts 24pt glue
        \subheading{#1}%
        \bigskip
        }
 
%QED box, from the TeXbook, p. 106.
%\redefine\qed{{\unskip\nobreak\hfil\penalty50\hskip2em\vadjust{}\nobreak\hfil
%    $\square$\parfillskip=0pt\finalhyphendemerits=0\par}}
 
%Other math symbols
        \define \x{\times}
        \let    \< = \langle
        \let    \> = \rangle
 
%Operator name macros
\define\op#1{\operatorname{\fam=0\tenrm{#1}}} %for text in math mode use 
                                              %\op{...}
\define\opwl#1{\operatornamewithlimits{\fam=0\tenrm{#1}}}
 
%Faye's AmS-TeX macros for this file
        %Greek letters
        \define         \a              {\alpha}
        \redefine       \b              {\beta}
        \redefine       \d              {\delta}
        \redefine       \D              {\Delta}
        \define         \e              {\varepsilon}
        \define         \g              {\gamma}
        \define         \G              {\Gamma}
        \redefine       \l              {\lambda}
        \redefine       \L              {\Lambda}
        \define         \n              {\nabla}
        \redefine       \var    {\varphi}
        \define         \s              {\sigma}
        \redefine       \Sig    {\Sigma}
        \redefine       \t              {\tau}
        \define         \th             {\theta}
        \redefine       \O              {\Omega}
        \redefine       \o              {\omega}
        \define         \z              {\zeta}
        
        %Other macros
        \redefine       \i              {\infty}
        \define         \p              {\partial}
 
% Fill in the blanks or delete the following
\topmatter
\title A Partition Theorem for Pairs of Finite Sets
\endtitle
\author Thomas Jech${}^1$ and Saharon Shelah${}^2$
\endauthor
\abstract{%
Every partition of $[[\o_1]^{<\o}]^2$ 
into finitely many pieces has a cofinal homogeneous set.  
Furthermore, it is consistent that every directed partially ordered set 
satisfies the partition property if and only if it has finite character.}
\endabstract
\address{%
T\. Jech, Department of Mathematics, The Pennsylvania State University, 
and Mathematical Sciences Research Institute, Berkeley, CA~~94720 \newline\newline 
S\. Shelah, Department of Mathematics, The Hebrew University, and 
Mathematical Sciences Research Institute, Berkeley, CA~~94720}\endaddress
\thanks{${}^1$Supported by NSF Grant DMS-8614447 and by BSF. \newline
${}^2$Supported at MSRI by NSF Grant 8505550 and by BSF. Publ\. No\. 392}
\endthanks
\endtopmatter
\baselineskip=14 pt
 
\document
 
\sectionskip

\newpage
 
\section{1.  Introduction}
 
A branch of combinatorics called {\it Ramsey theory} deals with phenomena
of the following kind: If a sufficiently 
large set of objects is endowed with a certain
structure then a large subset
can be found whose all elements are ``alike''.
 
A simple instance is the pigeon-hole principle: If there are more pigeons
than pigeon-holes then some pigeon-hole is occupied by more than one pigeon.
Another is this (the finite Ramsey theorem): For every integer $k>2$ there is
an integer $n$ with the property that if $A$ is a set of at least
$n$ elements and if the set of all (unordered) pairs $\{a,b\}\subset A$ is 
divided into two classes then there is a subset $H\subset A$ of size $k$
such that all pairs $\{a,b\}\subset H$ belong to the same class.
 
Many such principles have been formulated and proved, with applications
in various branches of mathematics; most are variants of Ramsey's Theorem
\cite{2}.
 
Ramsey's Theorem states 
(in particular) that every partition of the set $[\bold N]^2$ (into 
finitely many pieces) has an infinite homogeneous set, i\.e\. a set $H 
\subseteq 
\bold N$ cofinal in 
$(\bold N,<)$ such that $[H]^2$ is included in one piece of the 
partition.  The following generalization of Ramsey's Theorem was suggested in 
\cite{3}:
 
Let $A$ be an infinite set, and let $[A]^{<\o}$ denote the set of all finite 
subsets of $A$.  A set $H \subseteq [A]^{<\o}$ is {\it cofinal} in $[A]^{<\o}$ 
if for every $x \in [A]^{<\o}$ there exists a $y \in H$ such that $x \subseteq 
y$. Note that if a cofinal set $H$ is partitioned into two pieces,
$H = H_1 \cup H_2,$ then at least one of the two sets $H_1$, $H_2$ is cofinal.
 
Let $F: [[A]^{<\o}]^2 \to \{1,\dots,k\}$ be a partition of pairs of finite 
subsets of $A$; a set $H \subseteq [A]^{<\o}$ is {\it homogeneous} for $F$ if 
all pairs $(a,b) \in [H]^2$ with the property that $a \subset b$ belong to the 
same piece of the partition, i\.e\.
$$
F(x_1,x_2) = F(y_1,y_2)
$$
whenever $x_1,x_2,y_1,y_2 \in H$ and $x_1 \subset x_2$, $y_1 \subset y_2$.
 
The question raised in \cite{3} asked whether for every infinite $A$, every 
partition of $[[A]^{<\o}]^2$ has a cofinal homogeneous set.
 
It is not difficult to see that if $A$ is countable, then $[A]^{<\o}$ has a 
cofinal subset of order type $\o$ and so $[A]^{<\o}$ satisfies the partition 
property as a consequence of Ramsey's Theorem.  For an arbitrary $A$, the 
partition property in question is a generalization of Ramsey's Theorem for
pairs.
 
We answer the question in the affirmative in the case when $|A| = \aleph_1$:
 
\proclaim{Theorem 1} If $|A| = \aleph_1$, then every partition of 
$[[A]^{<\o}]^2$ into finitely many pieces has a cofinal homogeneous set.
\endproclaim
 
The question remains open for sets of size greater than $\aleph_1$.  By an 
unpublished theorem of Galvin, Martin's Axiom implies the partition property 
for all sets $A$ of cardinality less than $2^{\aleph_0}$.
 
More generally, let $S$ be a partially ordered set, and assume that $S$ is 
directed and does not have a maximal element.  
A set $H \subseteq S$ is {\it cofinal} in $S$ if for every $x \in S$ there
exists a $y \in H$ such that $x \le y.$ Let $r \ge 2$ and $k \ge 2,$ and
let $F: [S]^r \to \{1,\dots,k\}$ be a partition of $r$-tuples in $S.$
A set $H \subseteq S$ is {\it homogeneous} for $F$ if for all $x_1,\dots,
x_r$ and $y_1,\dots,y_r$ such that  $x_1 < \dots < x_r$ and
$y_1 < \dots < y_r$ we have
$$F(x_1,\dots,x_r) = F(y_1,\dots,y_r).$$
Using the standard arrow notation, the formula
$$S \rightarrow (\text{cofinal subset})^r_k$$
states that for every partition $F: [S]^r \to \{1,\dots,k\}$ there exists
a cofinal subset $H$ of $S$ homogeneous for F.
 
The following is an unpublished result of Galvin \cite{1}:
 
\proclaim{Theorem 2} (F. Galvin) Assume MA($\kappa$). Let $S$ be a partially
ordered set of power $\kappa$, which is directed, and suppose for all
$a \in S$, $\{b \in S : b < a\}$ is finite. Let $f : \{(x,y) \in S \x S :
x < y\} \to  \{\op{red, blue}\}.$ Then there is a cofinal $H \subseteq S$
such that $f$ is constant on $\{(x,y) \in H \x H : x < y\}.$
\endproclaim
 
Galvin's method admits a generalization to partitions of $r$-tuples, for any
$r \ge 2$ (see the proof of Theorem 4 below). Thus assuming Martin's Axiom
the following holds:
 
\proclaim{Theorem 2'} 
Let $S$ be a directed partially ordered set of cardinality less than 
$2^{\aleph_0},$ without maximal 
element and such that for every $a \in S$ the set $\{x \in S : x < a\}$
is finite. Then
$$S \rightarrow (\text{cofinal subset})^r_k \qquad \text{for all }r, k <\o.$$
\endproclaim
 
Note that every partially ordered set $S$ with the properties stated above 
is isomorphic to a cofinal subset of $[S]^{<\o}.$ 
 
The statement that for every cofinal $S \subseteq [\o_1]^{<\o},$
$$S \rightarrow (\text{cofinal subset})^2_2$$
is not a theorem of ZFC, as by an unpublished result of Laver \cite{4}
a counterexample exists under the assumption of the continuum hypothesis:
 
\proclaim{Theorem 3} (R. Laver) Let $\kappa$ be a cardinal such that
$\kappa^{\aleph_0} = \kappa.$ Then there exist a cofinal set $S \subset
[\kappa]^{<\o}$ and a partition $F: [S]^2 \to \{1,2\}$ such that no cofinal
subset of $S$ is homogeneous for $F.$
\endproclaim
 
With Laver's permission we include the proof of Theorem 3 below.
 
We say that a partially ordered set $S$ has {\it finite character} if $S$
has a cofinal set $S'$ such that every 
$x \in S'$ has only finitely many predecessors in $S'$.
Thus Galvin's theorem implies that
$$S \rightarrow (\text{cofinal subset})^r_k$$
holds for every set $S$ of size $\aleph_1$ that has finite character, if
Martin's Axiom holds together with $2^{\aleph_0} > \aleph_1,$ and Laver's
theorem implies that if $2^{\aleph_0} = \aleph_1$ then a partial order $S$
exists that has size $\aleph_1$ and finite character but
$$S \rightarrow (\text{cofinal subset})^2_2$$
fails.
 
\proclaim{Theorem 4} In the Cohen model for $2^{\aleph_0} = \aleph_2$ the
following statements are equivalent for every directed set of cardinality
$\aleph_1:$
\roster
\item $S \rightarrow (\text{cofinal subset})^2_2$
\item $S \rightarrow (\text{cofinal subset})^r_k$ for all $r, k < \o$
\item $S$ has finite character.
\endroster 
\endproclaim
 
The consistency proof of ``(3) implies (2)'' is essentially 
the Galvin's result; we will
show that (1) implies (3) in the Cohen model.

\section{2.  Proof of Theorem 1}
 
Throughout this section we consider a fixed partition $F:  
[[\o_1]^{<\o}]^2 \to \{1,\dots,k\}.$  The pairs $\{x,y\}$ 
such that $x \subset y$ are divided into two 
classes; we shall refer to these two classes as {\it colors}.
 
We reserve lower case letters such as $a$, $b$, $c$ for finite subsets of 
$\o_1$, and capital letters such as $A$, $B$, $C$ for at most
countable subsets of 
$\o_1$.
 
A {\it partial coloring} of a finite set $a$ is a function $f$ whose domain is 
a set of subsets of $a$, and whose values are in the set 
$\{1,\dots,k\}$.  A {\it total coloring} of $a$ is a partial 
coloring whose domain is the set of all subsets of $a$.
 
If $a \subset b$ and if $f$ is a partial coloring of $a$ then $b$ is $f$-{\it 
correct} if for every $x \in \op{dom}(f)$, the pair $(x,b)$ has the color 
$f(x)$ (i.e. $F(x,b)=f(x)$).
 
If $a \subseteq A$ then $b$ is an $A$-{\it extension} of $a$,
$a \le_A b$, if $a \subseteq b$ 
and $b \cap A = a$. An $A$-extension $b$ of $a$ is {\it proper}
if $a \subset b.$
 
We shall consider pairs $(a,A)$ where $a$ is finite, $A$ is at most countable
and $a \subseteq A.$
If $a \subseteq A$ and $b \subseteq B$ then
$$
(a,A) \le (b,B)
$$
means that $A \subseteq B$ and $b$ is an $A$-extension of $a$. Note that
$\le$ is transitive.
 
\demo{Definition 2.1} Let $a \subseteq A$, and let $f$ be a partial coloring 
of $a$.  We say that the pair $(a,A)$ is {\it good for} $f$ if for every 
$(b,B) \ge (a,A)$ there exists a 
proper $B$-extension $c$ of $b$ that is $f$-correct.
\enddemo
 
\demo{Remark} 
If $(a,A)$ is good for $f$ and if $f' \subseteq f$ and
$(a',A') \ge (a,A)$ then 
$(a',A')$ is good for $f'$.
\enddemo
 
\proclaim{Lemma 2.2} For every $(b,B)$ there exist a total coloring $g$
of $b$ and some $(c,C) \ge (b,B)$ such that $(c,C)$ is good for $g.$
 
Moreover, we may require that $c$ is a proper $B$-extension of $b$
and is $g$-correct.
\endproclaim
 
\demo{Proof} First assume that $g$ and $(c,C) \ge (b,B)$ are as claimed
in the first part of the lemma. Then there is some $d >_C c$ that is 
$g$-correct, and $(d, C \cup d)$ is good for $g.$ Hence it suffices
to find for each $(b,B)$ a total coloring $g$ of $b$ and some $(c,C)
\ge (b,B)$ good for $g.$
 
Thus assume that the lemma fails and let $(b,B)$ be such that for every
total coloring $g$ of $b$, no $(c,C) \ge (b,B)$ is good for $g.$     
 
There are finitely many total colorings $g_1,\dots,g_m$ of $b$.  
We construct a sequence $(b_i,B_i)$, $i = 1,\dots,m$ so that
$$
(b,B) \le (b_1,B_1) \le \dots \le (b_m,B_m)
$$
as follows:
 
As $ (b,B)$ is not good for $g_1$, there exists some 
$(b_1,B_1) \ge (b,B )$ such that no proper $B_1$-extension of $b_1$ is 
$g_1$-correct.   
 
Next, as  $(b_1,B_1)$ is not good for $g_2$, there 
exists some $(b_2,B_2) \ge (b_1,B_1 )$ such that no 
proper $B_2$-extension of $b_2$ 
is $g_2$-correct.
 
And so on.  For each $i = 1,\dots,m$, no proper $B_i$-extension of $b_i$ is 
$g_i$-correct.
 
Now let  $c$ be an arbitrary proper
$B_m$-extension of $b_m.$  Let us consider the 
following total coloring $g$ of $b$:
$$
g(x)\ =\ F(x,c) \quad \text{ (the color of $(x,c)$)}.
$$
We have  $g = g_i$ for some $i \le 
m$.  It is now clear that $c$ is a $g_i$-correct 
proper $B_i$-extension of $b_i$, a 
contradiction. \qed
\enddemo

\proclaim{Lemma 2.3} If $(a,A)$ is good for $f$, then for every $(b,B) \ge 
(a,A)$ there exists a total coloring $g$ of $b$ extending $f$, and some $(c,C) 
\ge (b,B)$ such that $c$ is a $g$-correct proper 
$B$-extension of $b$ and $(c,C)$ is good for $g$.
\endproclaim
 
\demo{Proof} The proof proceeds as in Lemma 2.2, the difference being that
we consider only the total colorings $g_1,\dots,g_m$
of $b$ that extend $f.$ 
After having constructed $(b_1,B_1) \le \dots \le (b_m,B_m),$ we find
(because $(a,A)$ is good for $f$ and $(a,A) \le (b_m,B_m))$ a proper
$B_m$-extension $c$ of $b_m$ that is $f$-correct. Then $g$ (defined as
above) extends $f$ and so $g = g_i$ for some $i \le m.$ The rest of 
the proof is as before. \qed
\enddemo
 
We shall use Lemma~2\.2 
and Lemma 2.3 to construct an end-homogeneous cofinal set $H \subseteq 
[\o_1]^{<\o}$.
 
\demo{Definition 2.4} A set $H$ is {\it end-homogeneous} if for all $x,y,z \in 
H$, if $x \subset y$ and $x \subset z$, then $(x,y)$ and $(x,z)$ have the same 
color.
\enddemo
 
Note that if $H$ is a cofinal end-homogeneous set, then one of the sets
$$
H_i = \{a \in H: F(a,x)=i \text{ for all $x \in H$ such that } a \subset x\}
\qquad (i=1,\dots,k)
$$
is cofinal, and is homogeneous. 
It follows that it suffices to construct a cofinal 
end-homogeneous set.
\medpagebreak
 
\demo{Definition 2.5} An {\it approximation} is a triple $(A,G,H)$ where $A$ 
is an infinite
countable subset of $\o_1$, $G$ and $H$ are disjoint cofinal subsets of 
$[A]^{<\o}$, $H$ is end-homogeneous, and for every $a \in G$, $(a,A)$ is good 
for $f^H_a$, where $f^H_a$ is the partial coloring of $a$ defined on $\{x \subset 
a: x \in H\}$ by
$$
f^H_a(x) = \text{the color of $(x,y)$, where $y$ 
is any $y \in H$ such that $x \subset y$}.
$$
\enddemo
 
Let
$$
(A,G,H) \le (A',G',H')
$$
mean that $A \subseteq A'$, $G \subseteq G'$ and $H \subseteq H'.$  
We want to construct an increasing sequence of approximations 
$(A_{\a},G_{\a},H_{\a})$, such that $\bigcup_{\a} A_{\a} = \o_1$.  Then $H = 
\bigcup_{\a} H_{\a}$ is an end-homogeneous set, cofinal in $[\o_1]^{<\o}$.
 
It is easy to verify that if $\l$ is a countable  limit ordinal, and if 
$(A_{\a},G_{\a},H_{\a})$, $\a < \l$, is an increasing sequence of 
approximations, then $(\bigcup_{\a} A_{\a}, \bigcup_{\a} G_{\a},
\bigcup_{\a} H_{\a})$ 
is an approximation.  Thus to complete the proof, it suffices to prove the 
following two lemmas:
 
\proclaim{Lemma 2.6} There exists an approximation. \endproclaim
 
\proclaim{Lemma 2.7} Let $(A,G,H)$ be an approximation and let $\xi \in 
\o_1 - A$ 
be arbitrary.  Then there is an approximation $({\bar A},{\bar G},{\bar H}) \ge 
(A,G,H)$ such that $\xi \in {\bar A}$.
\endproclaim
 
\demo{Proof of Lemma 2.6}
We construct $A$ as the union of a sequence $c_0 \subset c_1 
\subset \dots \subset c_n \subset \dots$ of finite sets, as follows.   
Let $b_0$ be an arbitrary finite subset of $\o_1.$ 
By Lemma 2.2. there exist a 
total coloring $g_0$ of $b_0$ and some $(c_0,C_0)$ 
such that $(c_0,C_0)$ is good for $g_0$
and $c_0 \supset b_0$ is $g_0$-correct.  
 
Now let $n \ge 0;$ we have constructed $(c_0,C_0),\dots,(c_n,C_n)$ such that
$c_0 \subset \dots \subset c_n .$ Fix for each $i\le n$ an enumeration
of $C_i$ of order-type $\o.$
Let $b_{n+1} \supseteq c_n$ be 
a finite set such that for each $i \le n,$ $b_{n+1}$ contains the first 
$n$ elements of  $C_i.$  
This will guarantee that $\bigcup_{n=0}^{\i} c_n = \bigcup_{n=0}^{\i} C_n$.
 
By Lemma 2.2. there exist a total coloring 
$g_{n+1}$ of $b_{n+1}$ and some $(c_{n+1},C_{n+1})$ 
such that $(c_{n+1},C_{n+1})$ is good for $g_{n+1}$ and $c_{n+1} \supset
b_{n+1}$ is $g_{n+1}$-correct. 
 
We let $A = \bigcup_{n=0}^{\i} C_n$.  To construct $G$ and $H$, 
consider  the partition $F$ restricted to the set 
$[\{c_n\}_{n=0}^{\i}]^2$.  By Ramsey's Theorem, $\{c_n\}_{n=0}^{\i}$ has an 
infinite homogeneous (let us say green) subsequence.  Let us denote this 
subsequence
$$
d_0 \subset e_0 \subset d_1 \subset e_1 \subset \dots \subset d_i \subset e_i 
\subset \dots
$$
and let $G =  \{d_i\}_{i=0}^{\i}$, $H =  
\{e_i\}_{i=0}^{\i}$.  Clearly, $G$ and $H$ are disjoint cofinal 
subsets of $[ A]^{<\o}$. Moreover, $ H$ is homogeneous, 
and we claim that for every $a \in G$, $(a, A)$ is good for $f^H_a$.
Since $a=c_n$ for some $n$, $(c_n,C_n)$ is good for $g_n$ and $(c_n,C_n)
\le (c_n,A),$ it suffices to show that $f^H_a \subseteq g_n.$
If $x \in \text{dom } f^H_a$ then because $g_n$ is a total coloring
of $b_n$ and $x=c_m$ for some $m < n,$ and because $c_n$ is
$g_n$-correct, we have $g_n(x) =$ the color of $(x,a),$ which is green,
because both $x$ and $a$ are in the homogeneous sequence.
But $f^H_a(x)$ is also green. Hence $f^H_a \subseteq g_n,$ and
$(A,G,H)$ is an approximation. \qed
\enddemo
 
\demo{Proof of Lemma 2.7} 
We construct $\bar A$ as the union of a sequence $c_0 \subset c_1 
\subset \dots \subset c_n \subset \dots$ of finite sets, as follows.  First, 
we choose an increasing cofinal sequence $a_0 \subset a_1 \subset \dots 
\subset a_n \subset \dots$ in $G$.  Let $b_0$ be some $A$-extension of $a_0$ 
such that $\xi \in b_0$.  As $(a_0,A)$ is good for $f^H_{a_0}$, there exist a 
total coloring $g_0$ of $b_0$ extending $f^H_{a_0}$, an $A$-extension $c_0$ of 
$a_0$ such that $c_0 \supset b_0$
and some $C_0 \supseteq A \cup c_0$ such that $c_0$ is $g_0$-correct and 
$(c_0,C_0)$ is good for $g_0$.  
 
Now assume that $(c_n,C_n)$ has been constructed and $c_n$ is an $A$-extension 
of $a_n$.  Let $b_{n+1}$ be some $A$-extension of $a_{n+1}$ such that $b_{n+1} 
\supseteq c_n$.  Moreover, we choose $b_{n+1}$ large enough to contain the first 
$n$ elements of each $C_i - A$, $i = 0,\dots,n$ (in some fixed enumeration).  
This will guarantee that $\bigcup_{n=0}^{\i} c_n = \bigcup_{n=0}^{\i} C_n$.
 
As $(a_{n+1},A)$ is good for $f^H_{a_{n+1}}$, there exist a total coloring 
$g_{n+1}$ of $b_{n+1}$ extending $f^H_{a_{n+1}}$, an $A$-extension $c_{n+1}$ of 
$a_{n+1}$ such that $c_{n+1} \supset b_{n+1}$
and some $C_{n+1} \supseteq A \cup c_{n+1}$ such that $c_{n+1}$ is 
$g_{n+1}$-correct and $(c_{n+1},C_{n+1})$ is good for $g_{n+1}$.
 
We let ${\bar A} = \bigcup_{n=0}^{\i} C_n$.  To construct ${\bar G}$ and ${\bar 
H}$, consider the partition restricted to the set 
$[\{c_n\}_{n=0}^{\i}]^2$.  By Ramsey's Theorem, $\{c_n\}_{n=0}^{\i}$ has an 
infinite homogeneous (let us say green) subsequence.  Let us denote this 
subsequence
$$
d_0 \subset e_0 \subset d_1 \subset e_1 \subset \dots \subset d_i \subset e_i 
\subset \dots
$$
and let ${\bar G} = G \cup \{d_i\}_{i=0}^{\i}$, ${\bar H} = H \cup 
\{e_i\}_{i=0}^{\i}$.  Clearly, ${\bar G}$ and ${\bar H}$ are disjoint cofinal 
subsets of $[{\bar A}]^{<\o}$, and $G = {\bar G} \cap [A]^{<\o}$, $H = {\bar 
H} \cap [A]^{<\o}$.  It remains to show that ${\bar H}$ is end-homogeneous, 
and that for every $a \in {\bar G}$, $(a,{\bar A})$ is good for $f^{\bar H}_a$.
 
To prove that ${\bar H}$ is end-homogeneous, we have to show that the color of 
$(x,y)$ for $x,y \in {\bar H}$ does not depend on $y$.  If $x \in {\bar H} - 
H$, say $x = e_i$, then every $y \supset x$ in ${\bar H}$ is some $e_m$, and 
$(e_i,e_m)$ is green.  If $x \in H$, then the color of $(x,y)$ is determined 
by $H$, and should be equal to $f^H_a(x)$, for any $a \supset x$ in $G$.  We 
have to show that $(x,e_i)$ has this color, for all $e_i \supset x$.  So let 
$i$ be such that $e_i \supset x$; we have $e_i = c_n$ for some $n$.  As $c_n$ 
is an $A$-extension of $a_n$, it follows that $x \subset a_n$.  Since 
$c_n$ is $g_n$-correct and $g_n \supseteq f^H_{a_n}$, $c_n$ is 
$f^H_{a_n}$-correct.  Therefore $(x,c_n)$ has color $f^H_{a_n}(x)$.
 
Finally, we prove that for every $a \in {\bar G}$, $(a,{\bar A})$ is good for 
$f^{\bar H}_a$.  If $a \in G$, then $f^{\bar H}_a$ is just  $f^H_a$ because $\{x \subset a: 
x \in {\bar H}\} = \{x \subset a: x \in H\}$.  Because $(a,A)$ is good for 
$f^H_a$, and $A \subseteq {\bar A}$, $(a,{\bar A})$ is good for 
$f^{\bar H}_a$.  So let 
$a \in {\bar G} - G$, say $a = d_i = c_n$.  We know that $(c_n,C_n)$ is good 
for $g_n$, and $C_n \subseteq {\bar A}$, so it suffices to show that 
$f^{\bar H}_a \subseteq g_n$, 
and then it follows that $(a,A)$ is good for $f^{\bar H}_a$.
 
So let $x \subset a$ be an element of ${\bar H}$.  If $x \in H$ then, because 
$a = c_n$ is an $A$-extension of $a_n$, $x \subset a_n$ and so $x \in 
\op{dom}(f^H_{a_n})$.  We already know that ${\bar H}$ is end-homogeneous, so 
$f^{\bar H}_a(x) = f^H_{a_n}(x) =$ 
the color of $(x,y)$ for any $y \supset x$ in ${\bar 
H}$.  Because $g_n$ is an extension of $f^H_{a_n}$, 
we have $f^{\bar H}_a(x) = g_n(x)$.
 
If $x \in {\bar H} - H$ then $x = c_m$ for some $m < n$, and 
$f^{\bar H}_a(x) =$ green (because $x=e_i$ for some $i$).  
Now $g_n$ is a total coloring of $b_n$, and $b_n \supseteq c_m$, so $x \in 
\op{dom}(g_n)$ and it remains to show that $g_n(x) =$ green.  But $c_n$ is 
$g_n$-correct, and so $g_n(x) =$ the color of $(c_m,c_n) =$ green. \qed
 
\enddemo
 
\section{3.  Proof of Theorem 4}
 
We shall prove that the equivalence between the partition property and finite 
character, for directed partial orders of size $\aleph_1$, holds in the model 
$V[G]$ obtained by adding $\aleph_2$ Cohen reals to a ground model 
for ZFC.
 
We shall first prove that (3) implies (1) in the Cohen model, and then outline
how the proof can be modified to show that (3) implies (2).
Assume that $S$ is a directed partially ordered set of size $\aleph_1$ 
in the model $V[G]$, and assume that each $a \in S$ has finitely many 
predecessors.  Let $F$ be a partition of $[S]^2$.  As 
$|S| = |F| = \aleph_1,$ $V[G]$ is a generic 
extension of $V[S,F]$ by the Cohen forcing, and we may 
assume that $S$ and $F$ are in the ground model.  Thus it suffices to prove 
that adding $\aleph_1$ Cohen reals produces a 
cofinal homogeneous set for $F$.
 
In fact, we define a forcing notion $P$ that produces a generic cofinal 
homogeneous set for $F$, and then show that $P$ is equivalent to adding 
$\aleph_1$ Cohen reals.  The forcing notion $P$ is essentially the one used by 
Galvin in his proof of the partition property for $[\o_1]^{<\o}$ from Martin's 
Axiom.
 
Let $D$ be an ultrafilter on $S$ with the property that for every $a \in S$, 
$\{x \in S: a \le x\} \in D$.  We say that $a \in S$ is {\it red}, if for 
$D$-almost all $x > a$, $(a,x)$ is red; otherwise $a$ is {\it green}.  
Either almost 
all $a \in S$ are red, or almost all are green; let us assume that almost all $a 
\in S$ are red.  A forcing condition in $P$ is a finite red-homogeneous set 
of red 
points.  A condition $p$ is stronger than $q$ if $p \supseteq q$ and if for no 
$x \in p - q$ and no $y \in q$ we have $x < y$.
 
Using the ultrafilter $D$ one can easily verify that for every $p \in P$ and 
every $a \in S$ there exists some $x \ge a$ such that $p \cup \{x\}$ is a 
condition stronger than $p$.  Therefore a generic set is a cofinal homogeneous 
set.
 
We shall finish the proof by showing that the forcing $P$ is equivalent to 
adding $\aleph_1$ Cohen reals.   
Let $S_{\a}$, $\a < \o_1$, be an elementary chain of countable 
submodels of $(S,<,\op{red})$, with limit $S$.  For each $\a$, let $P_{\a} = 
\{p \in P: p \subset S_{\a}\}$.  Each $P_{\a}$ is a countable forcing notion, 
therefore equivalent to adding a Cohen real.  It suffices to 
prove that every maximal antichain in $P_{\a}$ is a maximal antichain in $P$.  
This will follow from this claim:  For every $p \in P$ there is a ${\bar p} 
\in P_{\a}$ such that every $q \in P_{\a}$ stronger than ${\bar p}$ is 
compatible with $p$. Note that conditions $p$ and $q$ are compatible if
and only if no element of $p - q$ is less than any element of $q$ and no
element of $q - p$ is less than any element of $p$.
 
Let $p \in P$.  Let $Z$ be the (finite) set 
$\{x \in S_{\a}: x \le a$ for some $a \in p\}$, and let $u \in S_{\a}$ be a 
red point such that $u > x$ for all $x \in Z$, and that $(x,u)$ is 
red for all $x \in p \cap S_{\a}$.  Such a $u$ exists as $S_{\a}$ is an 
elementary submodel.  Now let ${\bar p} = (p \cap S_{\a}) \cup \{u\}$.
 
Clearly, ${\bar p}$ is a condition in $P_{\a}$.  Let $q \in P_{\a}$ be 
stronger than ${\bar p}$ and let us show that $q$ and $p$ are compatible.   
First, let $x \in q - p$ and $y \in p$. We claim that $x$ is not
less than $y$: since $q$ is stronger than ${\bar p}$,
$x$ is not less than $u$,
hence $x \notin Z$ and because $x \in S_{\a},$ the claim follows.
 
Second, let $x \in q$ and $y \in p -q$. We claim that $y$ is not
less than $x$: this is 
because $x \in S_{\a}$, $y \notin S_{\a}$, and since $S_{\a}$ is an elementary 
submodel and $x$ has finitely many predecessors, all $z < x$ are in $S_{\a}$. 
 
Hence $p$ and $q$ are compatible. 
\medpagebreak
 
We shall now outline how the above proof is modified to show that (3)
implies (2) in the Cohen model. For instance, let $k=2$ and $r=3$. The
above proof produces in fact a homogeneous cofinal set $H$ such
that $D \cup \{H\}$ has the finite intersection property. (For every
condition and every $A \in D$ there exists a stronger condition $q$ such
that $q \cap A \ne \emptyset.$)
 
Let $F$ be a partition of $[S]^3$ into \{red, green\}. For each $a \in S,$ 
let $F_a$ be the partition of $[S]^2$ given by $F_a(x,y) = F(a,x,y).$
Let $D$ be an ultrafilter on $S$ as before, and let $P_a$ denote the forcing
that produces a homogeneous cofinal set for $F_a.$ The product of
$\{P_a : a \in S\}$ is isomorphic to adding $\aleph_1$ Cohen reals and if
$\{H_a : a \in S\}$ are the generic homogeneous cofinal sets then
$D \cup \{H_a : a \in S\}$ has the finite intersection property.
 
We may therefore assume that the sets $H_a$ are in the ground model,
and $H_a \in D$ for each $a \in S.$ We say that $a \in S$ is {\it red},
if $H_a$ is red-homogeneous; otherwise $a$ is {\it green.} Assuming that 
almost all $a \in S$ are red, a forcing condition is a finite
red-homogeneous set of red points. This forcing produces a cofinal homogeneous
set for the partition $F$, and is equivalent to adding $\aleph_1$
Cohen reals. \qed
\medpagebreak
 
We shall now prove that (1) implies (3) in the Cohen model $V[G].$  So 
let $S \in V[G]$ be a directed partially ordered set of size $\aleph_1$ and 
assume that $S$ has the partition property.  Consider the forcing notion $P$ 
that adds, with finite conditions, a generic partition of $[S]^2$:
 
The forcing conditions in $P$ are functions whose domain is a finite
subset of $[S]^2$, with 
values \{red, green\}, and let ${\dot F}$ be the canonical name for a 
$P$-generic set.  Clearly, $P$ is equivalent to adding $\aleph_1$ Cohen reals, 
and if $Q$ is the forcing that adds $\aleph_2$ Cohen reals, we have $Q 
\x P \simeq Q$.  We shall prove:
 
\proclaim{Lemma} $P$ forces that if ${\dot F}$ has a cofinal homogeneous set, 
then $S$ has finite character.
\endproclaim
 
Granted the lemma, we complete the proof of Theorem~2 as follows:  Let ${\dot 
S}$ be a $Q$-name for $S \in V[G]$, and let $R$ be the forcing such that $V^Q 
= V[{\dot S}]^R$.  We have $R \simeq Q$ and so $R \simeq R \x P$.  The 
assumption is that $R$ (and therefore $R \x P$) forces that every partition of 
$S$ has a cofinal homogeneous set.  Hence $R \x P$ forces that ${\dot F}$ has 
a cofinal homogeneous set, and it follows from the lemma that $R \x P$ forces 
that $S$ has finite character.  Hence in $V[G]$, $S$ has finite character.
 
\demo{Proof of Lemma} Let ${\dot H}$ be a $P$-name for a cofinal homogeneous 
set for ${\dot F}$, and assume that $P$ forces that $[{\dot H}]^2$ is green.  
Let $S_{\a}$, $\a < \o_1$, be an elementary chain of countable submodels of 
$(S,<,P,\Vdash,{\dot F},{\dot Q})$.  First we claim that every condition 
forces the following:  For every $\a$, if $a \in {\dot H} - S_{\a}$ then the 
set $\{x \in {\dot H} \cap S_{\a}: x < a\}$ is finite.
 
So let us assume otherwise, and let $a \notin S_{\a}$ and $p \in P$ be such 
that $p \Vdash a \in {\dot H}$ and that $p \Vdash \{x \in {\dot H} \cap 
S_{\a}: x < a\}$ is infinite.  There is therefore some $x < a$, $x \in S_{\a}$ 
such that $(x,a) \notin \op{dom} p$ and that some $q$ stronger than $p$ forces 
$x \in {\dot H}$.  Since $S_{\a}$ is an elementary submodel, there 
is some $q$ stronger than 
the restriction of $p$ to $[S_{\a}]^2$ such that $\op{dom}(q) 
\subset [S_{\a}]^2$ and that $q$ forces $x \in {\dot H}$.  Now $q$ and $p$ are 
compatible conditions, and moreover, $(x,a)$ is not in the domain of $q \cup 
p$, so let $r$ be the extension of $p \cup q$ that forces that $(x,a)$ is red.  
Then $r \Vdash (x \in {\dot H}$ and $a \in {\dot H}$ and $(x,a)$ is red) which 
is a contradiction since $x < a$ and $[{\dot H}]^2$ is forced to be green. 
\enddemo
\medpagebreak
 
Now we shall construct, in $V^P$, a cofinal subset $C$ of $H$ such that each 
$a \in C$ has only finitely many predecessors in $C.$  
For each $\a$, let $a_{\a 0} 
\in S_{\a+1} - S_{\a}$ be, if it exists, an element of $H$ that is not below 
any $x \in H \cap S_{\a}$.  Then let $a_{\a n}$, $n < \o$, be an increasing 
sequence starting with $a_{\a 0}$, cofinal in $H \cap S_{\a+1}$.  Finally, let 
$C = \{a_{\a n}: \a < \o_1,\,n < \o\}$.
 
The set $C$ is cofinal in $H$.  If $a_{\a n} \in C$, then by the claim proved 
above, $a_{\a n}$ has only finitely many predecessors in $C \cap S_{\a}$, and 
because $a_{\b 0}$ is not less than $a_{\a n}$ for any $\b > \a$, $a_{\a n}$ has only 
finitely many predecessors in $C$. \qed
 
\section{4. Proof of Laver's Theorem}
 
Let $a_{\a}$ and $(M_{\a}, H_{\a})$, $\a < \kappa,$ enumerate, respectively,
the set $[\kappa]^{< \o}$ and the set of all pairs $(M, H)$ where
$M \in [\kappa]^{\le\aleph_0}$ and $H \subseteq [M]^{< \o}$ is cofinal
in $[M]^{< \o}.$ Furthermore, assume that $a_{\a} \subseteq \a$ and 
$M_{\a} \subseteq \a$ for all $\a.$
 
We construct a cofinal set $S = \{s_{\a} : \a < \kappa\}$ and a partition
$F : [S]^2 \to \{1,2\}$ as follows: Let $\a < \kappa.$ Let
$b_0 = a_{\a} \cup \{\a\};$ $\a$ is the largest element of $b_0.$
Choose, if possible, two distinct elements $c_0$ and $d_0$ of $H_{\a},$
and let $b_1 = b_0 \cup c_0 \cup d_0.$ Note that $\a$ is the largest
element of $b_1.$ Let $\a_1$ be the largest element of $b_1$ below
$\a,$ and choose, if possible, $c_1$ and $d_1$ in $H_{\a_1},$ distinct
from $c_0$ and $d_0$ and from each other, and let $b_2 = b_1 \cup
c_1 \cup d_1.$ Let $\a_2$ be largest in $b_2$ below $\a_1,$ and choose
$c_2$, $d_2$ in $H_{\a_2}$ distinct from $c_0$, $d_0$, $c_1$, $d_1$.
This procedure terminates after finitely many, say $k$, steps, and
we let $s_{\a} = b_k.$
 
For each $i \le k$, let $F(c_i, s_{\a}) = 1$ and $F(d_i, s_{\a}) =2$,
provided $c_i$ and $d_i$ are defined. Note that $\op{max} s_{\a} = \a$,
and that if $\b$ is the $i$ th  largest element of $s_{\a}$ and if $M_{\b}$
is infinite then $c_i$ and $d_i$ are defined; hence there exist $c$ and $d$
in $H_{\b}$ such that $F(c,s_{\a})=1$ and $F(d,s_{\a})=2.$
 
Let $S = \{s_{\a} : \a < \kappa\}$, and let $F$ be a partition  of $[S]^2$
that satisfies the conditions specified above. We claim that no cofinal
subset of $S$ is homogeneous for $F.$
 
Thus let $H$ be a cofinal subset of $S.$ There exists an infinite countable set
$M \subset \kappa$ such that $H \cap [M]^{< \o}$ is cofinal in $[M]^{< \o};$
let $\b < \kappa$ be such that $M_{\b} = M$ and $H_{\b} = H \cap [M]^{< \o}.$
As $H$ is cofinal, there is an $x \in H$ such that $\b \in x;$ as $H \subseteq
S$, there is some $\a$ such that $x = s_{\a}.$
 
Since $M_{\b}$ is infinite, there exist $c, d \in H_{\b}$ such that
$F(c, s_{\a})=1$ and $F(d,s_{\a})=2.$ Hence $H$ is not homogeneous for $F.$
\qed
 
\section{5. Open problems}
 
\roster
\item $[\aleph_2]^{< \o} \rightarrow (\text{cofinal subset})^2_2
\qquad$ (in ZFC)
 
\item $[\aleph_1]^{< \o} \rightarrow (\text{cofinal subset})^3_2
\qquad$ (in ZFC)
 
\item Is it consistent that there exists a directed partial ordering
of size $\aleph_1$ that does not have finite character but has the
partition property?
 
\item $[A]^{< \o} \rightarrow (\text{cofinal subset})^r_k$
\newline for all infinite sets $A$ and all integers $r, k \ge 2.$
\endroster
\heading References
\endheading
 
\vskip 0.5 truein
 
\ref \no 1
\by F\. Galvin
\paperinfo seminar notes from U.C.L.A
\endref
 
\ref \no 2
\by R\. Graham, B\. Rothschild and J\. Spencer
\book Ramsey Theory
\publ Wiley, New York
\yr 1980
\endref
 
\ref \no 3
\by T\. Jech
\paper Some combinatorial problems concerning uncountable cardinals
\jour Annals Math\. Logic \vol 5 \yr 1973 \pages 165--198
\endref
 
\ref \no 4
\by R\. Laver
\paperinfo private communication
\endref
 
\enddocument
\bye